\documentclass[a4paper, 12pt]{article}
\usepackage{adjustbox}
\usepackage{graphicx}
\usepackage{amsfonts}
\usepackage{amsbsy}
\usepackage{amssymb}
\usepackage[fleqn]{amsmath}
\usepackage{multicol}
\usepackage{longtable,ltxtable,booktabs}
\usepackage{blkarray}
\usepackage{afterpage}
\usepackage{float}
\usepackage{caption}
\usepackage{subcaption}
\usepackage{multirow}
\usepackage{pdflscape}
\usepackage{latexsym}
\usepackage{rotating}
\usepackage{enumerate}

\usepackage{mathtools}

\usepackage{epstopdf}
\usepackage{caption}
\usepackage{xcolor}
\usepackage[nodisplayskipstretch]{setspace}
\usepackage{amsmath}
\usepackage{longtable}
\usepackage[margin=.8 in]{geometry}

\usepackage{verbatim}

\usepackage{mathrsfs}
\usepackage{amsthm}
\usepackage{setspace}
 \usepackage{pdflscape}
\usepackage{colortbl}
\usepackage{latexsym}
\usepackage{amsfonts}
\usepackage{epsfig}
\usepackage{hyperref}
\usepackage{graphics}
\usepackage{float}
\theoremstyle{definition}
\newtheorem{thm}{Theorem}[section]

\newtheorem{defn}{Definition}[section]
\newtheorem{prop}{Proposition}[section]
\newtheorem{lemma}{Lemma}[section]

\theoremstyle{definition}

\title{ $p$-groups with small number of  character degrees and their normal subgroups}

\date{{}}\author{Nabajit Talukdar\thanks{Email address : ntalukdar2000@yahoo.co.in}\\
	\small and\\
	Kukil Kalpa Rajkhowa\thanks{Corresponding author. E-mail address : kukilrajkhowa@yahoo.com}\\
	\small Department of Mathematics\\
	\small Cotton University\\
	\small Guwahati-781001, India}

\begin{document}
\maketitle

\begin{abstract}
 If $G$ be a finite $p$-group and $\chi$ is a non-linear irreducible character of $G$, then $\chi(1)\leq |G/Z(G)|^{\frac{1}{2}}$. In \cite{fernandez2001groups}, Fern\'{a}ndez-Alcober and  Moret\'{o}  obtained the relation between the character degree set of a finite $p$-group $G$ and its normal subgroups depending on whether $|G/Z(G)|$ is a square or not. In this paper we  investigate the finite  $p$-group $G$ where for any normal subgroup $N$ of $G$ with $G'\not \leq N$ either $N\leq Z(G)$ or $|NZ(G)/Z(G)|\leq p$ and obtain some alternate characterizations of such groups. We find that if $G$ is a finite $p$-group with $|G/Z(G)|=p^{2n+1}$ and  $G$ satisfies the condition that for any normal subgroup $N$ of $G$ either $G'\not \leq N$ or $N\leq Z(G)$, then $cd(G)=\{1, p^{n}\}$. We also find that if $G$ is a finite $p$-group with nilpotency class not equal to $3$ and $|G/Z(G)|=p^{2n}$ and  $G$ satisfies the condition that for any normal subgroup $N$ of $G$ either $G'\not \leq N$ or $|NZ(G)/Z(G)|\leq p$, then $cd(G) \subseteq \{1, p^{n-1}, p^{n}\}$.
 
 \end{abstract}

  $\mathbf{2010\;Mathematics\;Subject\;Classification:}\;Primary: 20C15 $
	
	$\mathbf{Keywords\;and\;Phrases:}$ Character degrees, $p$-groups, nilpotency class

\section{Introduction}
In this paper, all groups  are finite. By $Irr(G)$ and $nl(G)$  we denote the set of of  complex irreducible characters and the set of   complex non-linear irreducible characters of the group $G$ respectively. By $c(G)$ we denote nilpotency class of $G$. The characterizations  of a finite group $G$ from the set $cd(G)$ of the degrees of its complex irreducible characters have been done in many papers. For example groups with two character degrees have been characterized in Chapter 12\cite{isaacs1994character} and Chapter 27 \cite{huppert1998character} . It can be found in Lemma 7.4\cite{huppert1998character} that for any $\chi \in Irr(G)$ we must get that $\chi(1)^{2}\leq |G/Z(G)|$. A group for which this bound is attained is called a group of central type. It has also been shown in \cite{isaacs1994character} that for any $\chi \in Irr(G)$ that $\chi(1)^{2}\leq |G/Z(\chi)|$ and if $G/Z(\chi)$ is Abelian then $\chi(1)^{2}=|G/Z(\chi)|$. The relation between the set of the degrees of its complex irreducible characters of a group and its normal subgroups have been studied in \cite{fernandez2001groups}.
We begin with the following definitions.
\begin{defn} \cite{fernandez2001groups}
A $p$-group is said to satisfy the strong condition on normal subgroups if for any $N\unlhd G$ either $G'\leq N$ or $N\leq Z(G)$.
\end{defn}

\begin{defn} \cite{fernandez2001groups}
A $p$-group is said to satisfy the weak condition on normal subgroups if for any $N\unlhd G$ either $G'\leq N$ or $|NZ(G)/Z(G)|\leq p$.
\end{defn}

%\begin{defn}
%\cite{nenciu2016nested}
%A non-Abelian group is called a generalized VZ-group (in short GVZ-group) if for all $\chi\in Irr(G)$ we get that $\chi(g)=0$ for all $g\in G \ Z(\chi)$.
%\end{defn}
%Remark/ By Corollary 2.30 of \cite{isaacs1994character} we get that in a GVZ-group $G$, $\chi(1)^{2}=|G/Z(\chi)|$ for all $\chi \in Irr(G)$.

In \cite{fernandez2001groups} the authors proved that for any $p$-group $G$, if $G$ satisfies the strong condition on normal subgroups then $c(G)\leq 3$ and if $G$ satisfies the weak condition on normal subgroups then $c(G)\leq 4$. Moreover in \cite{fernandez2001groups}, the character degree sets have been obtained for the $p$-groups satisfying the strong condition and weak condition based on  the value of $|G/Z(G)|$. In this paper we obtain some alternative characterizations of a $p$-group $G$ if $G$ satisfies the strong condition on normal subgroups and the weak condition on normal subgroups. We prove the following theorems.
%We prove the following results.\\

\begin{thm}
Let $G$ be a $p$-group of nilpotency class $2$ with two character degrees and  $|G/Z(G)|=p^{2n}$. Then the following conditions are equivalent:

   \begin{enumerate}[(i)]
        \item   $exp\ G'=p$ and $Z(\chi)/Z(G)$ is cyclic for all $\chi \in nl(G)$.

         %\item $cd(G)=\{1,p^{n}\}$.
         
         \item For every normal subgroup $N$ of $G$, either $G'\leq N$ or $N\leq Z(G)$.

         \item $(G,Z(G))$ is a generalized Camina pair.

       %  \item For every normal subgroup $N$ of $G$, either $G'\leq N$ or $|NZ(G)/Z(G)|\leq p$.
      \end{enumerate}
\end{thm}

\begin{thm}

Let $G$ be a finite $p$-group with $|G/Z(G)|=p^{2n+1}$. If $G$ satisfies the strong condition on normal subgroups then $cd(G)=\{1, p^{n}\}$. Moreover, if the nilpotency class of $G$ is $2$, then all the non-linear irreducible characters are non-faithful.
\end{thm}

%\begin{thm}

%Let $G$ be a finite $p$-group of nilpotency class $2$ and $|\frac{G}{Z(G)}|=p^{2n}$. Then the following conditions are equivalent/
%   \begin{enumerate}[(i)]

    %     \item $exp\ G'=p$, $|cd(G)|=2$ and $Z(\chi)/Z(G)$ is cyclic for all non linear irreducible character $\chi$ of $G$.
   
   %      \item $G$ satisfies the strong condition on normal subgroups.
  %       \item $G$ satisfies the weak condition on normal subgroups.
  %       \item $cd(G)=\{1, |G/Z(G)|^{\frac{1}{2}}\}$

 %   \end{enumerate}

%\end{thm}    

\begin{thm}
Let $G$ be a $p$-group of nilpotency class $2$ satisfying the weak condition on normal subgroups. Suppose one of the non-linear irreducible characters of $G$ is faithful. Then
   \begin{enumerate}[(i)]

     \item $G'=C_{p^{2}}$.

     \item  $K\cap Z(G) $ is the unique minimal normal $G$ for the kernel $K$ of any non-faithful non-linear irreducible character of $G$.

     \item $G$ has at least two distinct non-linear  non faithful irreducible  character kernels each of order  $p^{2}$.

      \item $K_{1}\cap K_{2}$ is the unique minimal normal subgroup of $G$ for any two distinct kernels $K_{1}$ and $K_{2}$ of  non-faithful non-linear irreducible characters of $G$.

      \item If $p=2$ then $G$ has exactly three distinct non-linear irreducible character kernels. 

    \end{enumerate}  
 \end{thm}

\begin{thm}

Let $G$ be a $p$-group of nilpotency class $2$ with two character degrees and $|G/Z(G)|=p^{2n+1}$. Then the following conditions are equivalent:

  \begin{enumerate}[(i)]
    
    \item   $exp\ G'=p$ and $Z(\chi)/Z(G)$ is cyclic for all non linear irreducible characters $\chi$ of $G$.

      %\item $cd(G)=\{1,p^{n}\}$.

     \item For every normal subgroup $N$ of $G$, either $G'\leq N$ or $|NZ(G)/Z(G)|\leq p$.

      \item $|Z(\chi)/Z(G)|=p$ for all non-linear irreducible characters $\chi$ of $G$.

     % \item $Z(G/N)=Z(\chi)/N$ for all normal subgroups $N$ with $G'\not \subseteq N$.

    \end{enumerate}  
\end{thm}

\begin{thm}

Let $G$ be a finite $p$-group of nilpotency class $2$ and $|G/Z(G)|=p^{2n}$. If $G$ satisfies the weak condition on normal subgroups then $cd(G) \subseteq \{1, p^{n-1}, p^{n}\}$.

\end{thm}

\section{Preliminaries}

%First we prove that if $G$ is a group satisfying the strong condition on normal subgroups then any non-Abelian factor group of $G$ also satisfies the strong condition on normal subgroups.
%\begin{lemma}
%\label{lemma1}
%Let $G$ be a group satisfying strong condition and $N$ be a normal subgroup of $G$ such that $G'\not \subseteq N$. Then the factor group $G/N$ also satisfies the strong condition.

%\begin{proof}
%  Let $L/N$ be a normal subgroup of $G/N$ such that $(G/N)'=G'N/N\not \subseteq L/N$. Then $G'\not \subseteq L$. By the stong condition in $G$ we get  that $L\leq Z(G)$. This gives that $L/N\leq Z(G)/N \leq Z(G/N)$.
%\end{proof}
  
%\end{lemma}

%First we solve a problem set by I.M. Isaacs in \cite{isaacs1994character}.
In this section, we prove some preliminary results which will aid us in proving our main results. First we prove some results for a $p$-group $G$ if $|G'|=p$. These results are stated in \cite{isaacs1994character} and \cite{doostie2012finite}.
\begin{lemma}
\label{lemma1}
Let $G$ be a $p$-group such that $|G'|=p$. Then 
\begin{enumerate}[(i)]
\item cd$(G)=\{1, |G/Z(G)|^{\frac{1}{2}}\}$.
\item $G$ satisfies the strong condition on normal subgroups.
\end{enumerate}
\begin{proof}
\begin{enumerate}[(i)]
 \item

Let $\chi$ be a non linear irreducible character and let $\Psi$ be the representation affording the character $\chi$.
We know $Z(G)\subseteq Z(\chi)$.
If possible let $g\in Z(\chi)-Z(G)$.
Then there exists $ h\in G$ such that $gh\neq hg$ that is $[g,h](\neq 1)\in G'$.
Since $|G'|=p$, $G'$ is cyclic and $G'=<[g,h]>$.
We have that $\Psi(g)\Psi(h)=\Psi(h)\Psi(g)$ that is $[g,h]\in ker\Psi=ker\chi$ and hence $G'=<[g,h]> \leq ker\chi$.
This gives that $\chi$ is a linear character and this contradiction proves that $Z(G)=Z(\chi)$.
Since $G'\subseteq Z(G)=Z(\chi)$, $G/Z(\chi)$ is Abelian and then by Theorem 2.31 in \cite{isaacs1994character} we get that $\chi(1)^{2}=|G/Z(\chi)|=|G/Z(G)|$.

\item 
  Let $N$ be a normal subgroup of $G$ such that $G'\not \subseteq N$. Then $N\cap G'=1$. For any $g\in G$ and $n\in N$, $g^{-1}n^{-1}gn\in N\cap G'=1$. Thus $gn=ng$ for all $n\in N, g\in G$. Hence $N\leq Z(G)$.
 \end{enumerate}

\end{proof}

\end{lemma}

From Theorem 2.32 in \cite{isaacs1994character} we get that if a group $G$ has a faithful irreducible character, then $Z(G)$ is cyclic. The necessary and sufficient conditions when all the non-linear irreducible characters are faithful have been obtained by Doostie and Saeidi in \cite{doostie2012finite}. We state the result here.

\begin{lemma}
 \cite{doostie2012finite}
\label{lemma2}

All non-linear irreducible characters of a finite group $G$ are faithful if and only if $|G'|=p$ and $Z(G)$ is cyclic.\\

\end{lemma}

Combining the above two lemmas we get the following result.

\begin{prop}

Let $G$ be a $p$-group where all the non-linear irreducible characters are faithful. Then cd$(G)=\{1, |G/Z(G)|^{\frac{1}{2}}\}$.

\end{prop}

If $G$ be a group of nilpotency class $2$ and if $G$ satisfies either the strong condition or the weak condition on normal subgroups then the exponents of the subgroup $G'$ and the factor group $G/Z(G)$ have been obtained in \cite{fernandez2001groups}.

\begin{thm}
\cite{fernandez2001groups}
\label{thm1}
Let $G$ be  $p$-group of nilpotency class $2$.
   \begin{enumerate}[(i)]
      \item If $G$ satisfies the strong condition on normal subgroups then $exp \ G/Z(G)=exp\ G'=p$.
      \item $G$ satisfies the weak condition on normal subgroups then $exp \ G/Z(G)=exp\  G'=p$ or $p^{2}$. In the later case $G/Z(G)\cong C_{p^{2}} \times C_{p^{2}}$ and $G'\cong C_{p^{2}}$.
   \end{enumerate}
\end{thm}

   In \cite{wang2008note} Qian and Wang examined the conditions under which the kernels of the non-linear characters of a non Abelian $p$-group form a chain with respect to inclusion and proved the following result. By $Kern(G)$ we denote the set of kernels of non-linear irreducible characters of $G$.

\begin{thm}

\label{thm2}
 Let $G$ be a finite non-Abelian $p$-group. Then the following are statements are equivalent:

 \begin{enumerate}[(i)]
     \item  $Kern(G)$ is a chain with respect to inclusion
      \item Whenever $N< G'$ is a normal subgroup of $G$, $N$ is a member of $Kern(G)$.
     \item $G$ is one of the following groups:
             \begin{enumerate}
                   \item  $G'$ is a unique minimal normal subgroup of $G$.
                   \item $G$ is of maximal class.

                   \end{enumerate}
    \end{enumerate} 

\end{thm}  

Next we prove that a $p$-group  of nilpotency class $2$ satisfying the strong condition must have that all the non linear  irreducible characters are either faithful or non-faithful.

\begin{lemma}

\label{lemma3}

Let $G$ a $p$-group of nilpotency class $2$ satisfying the strong condition. Then the non linear  irreducible characters are either  all faithful or all non-faithful.

\begin{proof}

 Suppose $G$ is $p$-group of nilpotency class $2$ satisfying the strong condition. Suppose at least one of the  non-linear irreducible characters of $G$ is faithful.  Then $Z(G)$ is cyclic and since the nilpotency class of $G$ is $2$, we get that $G'\subseteq Z(G)$. Since $Z(G)$ is a cyclic $p$-group, its subgroups form a chain with respect to inclusion. Since the group satisfies the strong condition, we get that the kernels of the irreducible characters are subsets of $Z(G)$ and hence it follows that all the non linear character kernels form a chain with respect to inclusion. Thus by Theorem \ref{thm2} either  $ G$ is of maximal class  or $G'$ is a unique minimal normal subgroup. If $G$ is of maximal class , then $|G|=p^{3}$ and hence all the character kernels are faithful. Suppose $G'$ is a unique minimal normal subgroup. Since $Z(G)$  is cyclic and $G'\subseteq Z(G)$ we get that $|G'|=p$. Then by Lemma \ref{lemma2} all the non linear character kernels are faithful.

\end{proof}

\end{lemma}

It has been proved in \cite{fernandez2001groups} that if the group $G$ satisfies the strong condition on normal subgroups then $c(G)\leq 3$ and  if $G$ satisfies the weak condition on normal subgroups then $c(G)\leq 4$. Here we provide new proofs of these results.

\begin{lemma}

\label{lemma3a}

Suppose $G$ be a non Abelian $p$-group which satisfies the strong condition on normal subgroups. Then $c(G)\leq 3$.

\begin{proof}

Suppose $c(G)>2$. Then $Z(G)< Z_{2}(G)$ and consequently by the strong condition we get that $G'\subseteq Z_{2}(G)$. Thus $G/Z_{2}(G)$ is Abelian and hence $c(G)=3$.

\end{proof}

\end{lemma}

\begin{lemma}

\label{lemma4}
Suppose $G$ be a non Abelian $p$-group which satisfies the weak condition on normal subgroups. Then $c(G)\leq 4$.

\begin{proof} 
Suppose $c(G)>2$. Then $Z(G)< Z_{2}(G)$.
First we assume that $G' \subseteq Z_{2}(G)$. Then $G/Z_{2}(G)$ is Abelian. Hence $c(G)=3$. Next we consider the case that  $Z_{2}(G) \not \subseteq G'$. We get that $|Z_{2}(G)/ Z(G)|\leq p$. Since $|Z_{2}(G)/ Z(G)|\neq 1$, we get that $|Z_{2}(G)/ Z(G)|= p$.
Now we must get that $G'\subseteq Z_{3}(G)$ because otherwise by the weak condition it follows that  $|Z_{3}(G)/ Z(G)|=p$ and hence $Z_{3}(G)=Z_{2}(G)$,  a contradiction. Thus $G/Z_{3}(G)$ is Abelian and consequently we get that $c(G)=4$.

\end{proof}

\end{lemma}

%%\begin{lemma}
%\label{lemma5}
%Let $G$ be a non-Abelian $p$ group which satisfies the weak condition on normal subgroups. Then the group $G/Z(G)$ satisfies the strong condition on normal subgroups.

%\begin{proof}

%Suppose $N/Z(G)\unlhd G/Z(G)$ such that  $(G/Z(G))'=G'Z(G)/Z(G) \not \subseteq N/Z(G)$. Then $G'\not \subseteq N$. By the weak condition we get %$|N/Z(G)|\leq p$.
%Since $G$ is a nilpotent group, if $|N/Z(G)|= p$ we get that $[N/Z(G), G/Z(G)]< N/Z(G)$.
%This gives that $[N/Z(G), G/Z(G)]=1$  and hence $N/Z(G)\subseteq Z(G/Z(G))$.

%\end{proof}

%\end{lemma}

A pair $(G,N)$ is said to be a generalized Camina pair (abbreviated GCP) if $N$ is normal in $G$ and all the non-linear irreducible characters of $G$ vanish outside $N$. 
The notion of GCP was introduced by Lewis in \cite{lewis2009vanishing}. An equivalent condition for a pair $(G,N)$ to be a GCP is: A pair $(G,N)$ is a GCP if and only if for $g\in G\setminus N$, the conjugacy class of $g$ in $G$ is $gG'$.
In \cite{prajapati2017irreducible} the authors obtained the following results.

\begin{thm}

\label{thm3}
Let $(G,Z(G))$ be a GCP. Then we have the following.

   \begin{enumerate}[(1)]

       \item $cd(G)=\{1, |G/Z(G)|^{\frac{1}{2}}\}$.

       \item The number of non-linear irreducible characters of $G$ is $|Z(G)|-|Z(G)/G'|$.

   \end{enumerate}
\end{thm}

\section{Main Results}

In \cite{fernandez2001groups} the authors characterized the $p$-groups $G$ satisfying the strong condition on normal subgroups by its character degree set if $G/Z(G)$ is an even power of $p$ . In the following theorem we provide some alternate characterizations. 

 \begin{thm}
Let $G$ be a $p$-group of nilpotency class $2$ with two character degrees and $|Z/Z(G)|=p^{2n}$. Then the following conditions are equivalent: 

    \begin{enumerate}[(i)]

        \item   $exp\ G'=p$ and $Z(\chi)/Z(G)$ is cyclic for all $\chi \in nl(G)$.

         %\item $cd(G)=\{1,p^{n}\}$.
         
         \item For every normal subgroup $N$ of $G$, either $G'\leq N$ or $N\leq Z(G)$.

         \item $(G,Z(G))$ is a generalized Camina pair.

     \end{enumerate}

  \begin{proof} $(i)\Rightarrow (ii)$ 
    Since the nilpotency class of $G$ is $2$ and exp $G'=p$, it follows from Lemma 4.4 \cite{isaacs2008finite} that exp $(G/Z(G))=p$ and hence $G/Z(G)$ is elementary Abelian. If $|G'|=p$, it follows from Lemma \ref{lemma1} that  $cd(G)=\{1, |G/Z(G)|^\frac{1}{2}\}$. So we assume that $|G'|>p$. Let us choose a normal subgroup $N$ of $G$ such that $[G'/N]=p$.  Let us consider an non linear irreducible  character $\chi \in Irr(G/N)$. Then $N\subseteq Ker\ \chi=K$ and hence $N=K\cap G'$. Hence $|(G/K)'|=[G'/ K\cap G']=p$. By lemma 2.2 we get 
    cd$(G/K)=\{1, |G/Z(\chi)|^\frac{1}{2}\}$. Combining the fact that  $|G/Z(G)|=p^{2n}$ and $|G/Z(\chi)|$ has square order we get that $|Z(\chi)/Z(G)|$ has square order. Since exp$(Z(\chi)/Z(G)=p$ and $Z(\chi)/Z(G)$ is cyclic, we get that $Z(\chi)=Z(G)$. Hence $cd(G)=cd(G/K)=\{1, |G/Z(G)|^\frac{1}{2}\}$. Thus we get that $cd(G)=\{1, |G/Z(G)|^\frac{1}{2}\}$ and hence the result follows from Theorem B in 
 \cite{fernandez2001groups}.\\

     %Let $K$ be the kernel of a non-linear irreducible character of $G$. First assume that $|KZ(G)/Z(G)|= p$ . Since $KZ(G)/Z(G)\leq Z(\chi)/Z(G)$, by lemma we get that the Abelian group $Z(\chi)/Z(G)$ contains a subgroup isomorphic to $Z(\chi)/KZ(G)$, which is cyclic. Since $|KZ(G)/Z(G)|= p$, either $|Z(\chi)/KZ(G)|=1$  or $Z(\chi)/KZ(G)$  is of maximal order. If $|Z(\chi)/KZ(G)|=1$, then we get $Z(\chi)=Z(G)$.  If $|Z(\chi)/KZ(G)|$  is of maximal order then by the proof of 2.1.2 of \cite{kurzweil2004finite} we get that $Z(\chi)/Z(G)\cong C_{p^{k}}\times C_{p}$, where $p^{k}=|Z(\chi)/KZ(G)|$. Since $G$ satisfies the weak condition on normal subgroups  by Theorem  we get that $exp\ G'=exp\ G/Z(G)=p$ or $p^{2}$. This gives that $exp\ Z(\chi)/Z(G)=p$ or $p^{2}$. If $exp\ Z(\chi)/Z(G)=p$ we get that $k=0$ or $k=1$.If $k=1$, then $|Z(\chi)/Z(G)|=p^{2}$ and hence $|G/Z(\chi)|=p^{2n-2}$. This contradicts that $|G/Z(\chi)|=\chi(1)^{2}$. Thus $k=0$ and hence 
  %$Z(\chi)/Z(G)\cong  C_{p}$. Since $|KZ(G)/Z(G)|=p$ and $KZ(G)/Z(G) \leq Z(\chi)/Z(G)$, we get that $KZ(G)/Z(G)=Z(\chi)/Z(G)$ and this in turn gives that $|Z(\chi)/Z(G)|=p$. 
  
 % Next assume that $|KZ(G)/Z(G)|= 1$ i.e. $K\leq Z(G)$. Then $Z(\chi)/Z(G)$ is cyclic and $exp\ Z(\chi)/Z(G)=p$. This gives that either $Z(\chi)=Z(G)$ or $|Z(\chi)/ Z(G)|=p$. If $Z(\chi)=Z(G)$ we get that $|\frac{G}{Z(\chi)}|=|\frac{G}{Z(G)}|=p^{2n+1}$, which contradicts that $(\chi(1))^{2}=|\frac{G}{Z(\chi)}|$. Hence $|Z(\chi)/Z(G)|=p$.

 $(ii)\Rightarrow (i)$ By Theorem \ref{thm1} we get that exp $G'=p$. Let $\chi \in nl(G)$ and  $K$ be the kernel of $\chi$. Then $Z(\chi)/K$ is a cyclic group. Since $K\leq Z(G)$ we get that $Z(G)/K$ is also a cyclic group. Hence it follows that $Z(\chi)/Z(G)$ is a cyclic group.\\

 $(i)\Rightarrow (iii)$ We get that $cd(G)=\{1, |G/Z(G)|^\frac{1}{2}\}$. Let $\chi \in nl(G)$. Since the group $G/Z(\chi)$ is Abelian, it follows by Theorem 2.31 in \cite{isaacs1994character} that $\chi(1)^{2}=|G/Z(\chi)|=|Z/Z(G)|$. Hence we get that $Z(\chi)=Z(G)$ . Hence by Corollary 2.30 in \cite{isaacs1994character} it follows that $\chi$ vanishes off $Z(G)$. Thus $(G,Z(G))$ is a generlized Camina pair.\\

 $(iii)\Rightarrow (ii)$  It follows from Theorem \ref{thm3} that $cd(G)=\{1,p^{n}\}$ and hence by theorem B of \cite{fernandez2001groups} we get that for any normal subgroup $N$ of $G$, either $G'\leq N$ or $N\leq Z(G)$.

\end{proof}  

 \end{thm}

In the following theorem we obtain the character degree set of a $p$ group $G$  under the condition that  $|G/Z(G)|=p^{2n+1}$ if $G$ satisfies the strong condition on normal subgroups.

\begin{thm}

Let $G$ be a finite $p$-group with $|G/Z(G)|=p^{2n+1}$. If $G$ satisfies the strong condition on normal subgroups then $cd(G)=\{1, p^{n}\}$. Moreover, if the nilpotency class of $G$ is $2$, then all the non-linear irreducible characters are non-faithful.

\begin{proof}
By Lemma \ref{lemma3a} we get that $c(G)\leq 3$.\\

Case I: nilpotency class of $G$ is $2$.\\

We get that $Z(\chi)/K$ is a cyclic group and since $K\leq Z(G)$ the subgroup $Z(G)/K$ is also cyclic. Hence $Z(\chi)/Z(G)$ is a cyclic group. Since $|G/Z(G)|=p^{2n+1}$ and $\chi(1)^{2}=|G/Z(\chi)|$, we get that $|Z(\chi)/Z(G)|$ is an odd power of $p$. By Theorem 4.3 of 
\cite{fernandez2001groups} we get that $exp\ Z(\chi)/Z(G)=p$. Hence 
$|Z(\chi)/Z(G)|=p$ for all non-linear irreducible characters of $G$. This gives that all the non-linear irreducible characters are non-faithful beacuse if $\chi$ is a faithful non-linear irreducible character of $G$ then by Lemma 2.27 of \cite{isaacs1994character} we get that $Z(\chi)=Z(G)$.\\
Since $\chi(1)^{2}=|G/Z(\chi)|$, it follows that $cd(G)=\{1, p^{n}\}$.

Case II: nilpotency class of $G$ is $3$.\\

By Theorem F of \cite{fernandez2001groups} we get that $|G/Z(G)|=p^{3}$ and hence $cd(G)=\{1,p\}$.

\end{proof}

\end{thm}

The structures of Abelian groups have been discussed in \cite{rotman2012introduction}. In particular we obtain the following result. 

\begin{lemma}
\label{lemma3.1}
Let $G$ be a finite Abelian group and $H\leq G$. Then $G$ contains a subgroup isopmorphic to $G/H$.
\end{lemma}

\begin{thm}

\label{thm3.3}
Let $G$ be a $p$-group of nilpotency class $2$ satisfying the weak condition on normal subgroups. Suppose one of the non-linear irreducible characters of $G$ is faithful. Then
   \begin{enumerate}[(i)]

     \item $G'=C_{p^{2}}$.

     \item  $K\cap Z(G) $ is the unique minimal normal $G$ for the kernel $K$ of any non-faithful non-linear irreducible character of $G$.

     \item $G$ has at least two distinct non-linear  non faithful irreducible  character kernels each of order  $p^{2}$.

      \item $K_{1}\cap K_{2}$ is the unique minimal normal subgroup of $G$ for any two distinct kernels $K_{1}$ and $K_{2}$ of  non-faithful non-linear irreducible characters of $G$.

      \item If $p=2$ then $G$ has exactly three distinct non-linear irreducible character kernels. 

    \end{enumerate}  
\begin{proof}

  Since one of the non-linear irreducible characters is faithful, we get that $Z(G)$ is cyclic.

\begin{enumerate}[(i)]

     \item By Theorem D of \cite{fernandez2001groups} we get that either $exp\ G'=p$ or $G'=C_{p^{2}}$. Since the nilpotency class of $G$ is $2$ and $G$ has a faithful non-linear irreducible character, we get that $Z(G)$ is cyclic. If $exp\ G'=p$ then $G'=C_{p}$.  
     
      Let $K$ be the kernel of a non-faithful non-linear irreducible character $\chi$. Now if $G'\subseteq K\cap Z(G)$ we get that $\chi$ is a faithful character . Hence we must get that $K\cap Z(G) < G'$. If  $G'=C_{p}$, this will give that $K\cap Z(G)=1$, a contradiction. Thus it follows  that $G'=C_{p^{2}}$.

    \item As observed in (i) above we get that $K\cap Z(G)< G'$ for the kernel $K$ of any non-faithful non-linear irreducible character of $G$ and hence it follows that $| K\cap Z(G) | =p$. Thus $K\cap Z(G)$ is the unique minimal normal subgroup of $G$.

      \item  
      Let $N=K\cap Z(G)$. Then we have that $|(G/N)'|=p$ . Hence by Lemma \ref{lemma1} we get that $G/N$ is a group of nilpotency class $2$ satisfying the strong condition on normal subgroups. Hence all the non-linear irreducible characters of $G/N$ are either faithful or non-faithful. If all the non-linear irreducible characters of $G/N$ are faithful we get that $N=K$ for the kernel $K$ of each of the non-faithful non-linear irreducible  characters of $G$. Thus $G$ is a $p$-group of nilpotency class $2$ where the number of distinct non-linear irreducible character kernels is $2$. Hence by Theorem 1.1 of \cite{doostie2012finite} we get that $G$ is a group of order $p^{4}$ and nilpotency class $3$, a contradiction. Thus all the non-linear irreducible characters of $G/N$ are non faithful and hence $N=K\cap Z(G) \subset K$ for the kernel $K$  of any non-faithful non-linear irreducible character of $G$. This gives that $|K|=p^{2}$.
      If $G$ has only one distinct kernel $K$ of non-faithful non-linear irreducible character, then by Theorem 1.1 and Lemma 2.2 of \cite{doostie2012finite} all the non-linear irreducible characters of $G/N$  are faithful. This contradiction proves that $G$ has at least two distinct non-linear  non faithful irreducible  character kernels.

       \item  Since $N \subseteq K_{1}\cap K_{2}$  and  $|N|=|K_{1}\cap K_{2}|=p$ we get that $K_{1}\cap K_{2}$ is the unique minimal normal subgroup of $G$ for any two distinct kernels $K_{1}$ and $K_{2}$ of  non-faithful non-linear irreducible characters of $G$.

       \item The result follows from Proposition 3.2 of \cite{doostie2012finite}.
       
\end{enumerate}
\end{proof}
\end{thm}

  In the following theorem we obtain some alternative characterizations of $p$ groups of nilpotency class 2 which satisfies the weak condition.

\begin{thm}

Let $G$ be a $p$-group of nilpotency class $2$ with two character degrees and $|G/Z(G)|=p^{2n+1}$. Then the following conditions are equivalent:

  \begin{enumerate}[(i)]
    
    \item   $exp\ G'=p$ and $Z(\chi)/Z(G)$ is cyclic for all non linear irreducible characters $\chi$ of $G$.

      %\item $cd(G)=\{1,p^{n}\}$.

     \item For every normal subgroup $N$ of $G$, either $G'\leq N$ or $|NZ(G)/Z(G)|\leq p$.

      \item $|Z(\chi)/Z(G)|=p$ for all non-linear irreducible characters $\chi$ of $G$.

     % \item $Z(G/N)=Z(\chi)/N$ for all normal subgroups $N$ with $G'\not \subseteq N$.

    \end{enumerate} 
\begin{proof}

   $(i)\Rightarrow (ii)$

   Since the nilpotency class of $G$ is $2$ and exp$G'=p$, it follows that exp$(G/Z(G))=p$ and hence $G/Z(G)$ is elementary Abelian. If $|G'|=p$, it follows from Lemma \ref{lemma1}  that  $\chi(1)^{2}=|G/Z(G)|$. This contradicts that $|G/Z(G)|=p^{2n+1}$.  So we get that $|G'|>p$. Let us choose a normal subgroup $N$ of $G'$ such that $[G'/N]=p$. Then $G/N$ is a non Abelian group and $cd(G/N)=cd(G)$. Let us consider an non linear irreducible  character $\chi \in Irr(G/N)$. Then $N\subseteq Ker\ \chi=K$ and hence $N=K\cap G'$. Hence $|(G/K)'|=[G'/ K\cap G']=p$. By Lemma \ref{lemma1} we get 
    cd$(G/K)=\{1, |G/Z(\chi)|^\frac{1}{2}\}$. Combining the fact that  $|G/Z(G)|=p^{2n+1}$ and $|G/Z(\chi)|$ has square order we get that $|Z(\chi)/Z(G)|$ is an odd power of $p$. Since exp $(Z(\chi)/Z(G)=p$ and $Z(\chi)/Z(G)$ is cyclic, we get that $|Z(\chi)/Z(G)|=p$. Hence $cd(G)=cd(G/K)=\{1, |G/Z(\chi)|^\frac{1}{2}\}=\{1,p^{n}\}$. Hence by Theorem C of \cite{fernandez2001groups} we get that for every normal subgroup $N$ of $G$, either $G'\leq N$ or $|NZ(G)/Z(G)|\leq p$. Thus the result follows from Theorem C in \cite{fernandez2001groups}.\\

  $(ii)\Rightarrow (iii)$ 
  Let $K$ be the kernel of a non-linear irreducible character of $G$.
  First assume that $|KZ(G)/Z(G)|= p$ for the kernel $K$ of a non-linear irreducible character of $G$. Since $KZ(G)/Z(G)\leq Z(\chi)/Z(G)$, by Lemma \ref{lemma3.1} we get that the Abelian group $Z(\chi)/Z(G)$ contains a subgroup isomorphic to $Z(\chi)/KZ(G)$, which is cyclic. Since $|KZ(G)/Z(G)|\leq p$, either $|Z(\chi)/KZ(G)|=1$  or $Z(\chi)/KZ(G)$  is of maximal order. If $|Z(\chi)/KZ(G)|=1$, then we get that $Z(\chi)=KZ(G)$ and hence $Z(\chi)/Z(G)|=|KZ(G)/Z(G)|=p$.  If $|Z(\chi)/KZ(G)|$  is of maximal order then by the proof of 2.1.2 of \cite{kurzweil2004finite} we get that $Z(\chi)/Z(G)\cong C_{p^{k}}\times C_{p}$, where $p^{k}=|Z(\chi)/KZ(G)|$. Since $|G/Z(G)|=p^{2n+1}$ and $G$ satisfies the weak condition on normal subgroups, by Theorem \ref{thm1} we get that $exp\ G'=exp\ G/Z(G)=p$. This gives that $exp\ Z(\chi)/Z(G)=p$ and consequently $k=0$ or $k=1$. If $k=1$, then $|Z(\chi)/Z(G)|=p^{2}$ and hence $|G/Z(\chi)|=p^{2n-1}$. This contradicts that $|G/Z(\chi)|=\chi(1)^{2}$. Thus $k=0$ and hence 
  $Z(\chi)/Z(G)\cong  C_{p}$. Since $|KZ(G)/Z(G)|=p$ and $KZ(G)/Z(G) \leq Z(\chi)/Z(G)$, we get that $KZ(G)/Z(G)=Z(\chi)/Z(G)$ and this in turn gives that $|Z(\chi)/Z(G)|=p$. Next assume that $|KZ(G)/Z(G)|= 1$ so that $K\leq Z(G)$. Then $Z(\chi)/Z(G)$ is cyclic and exp $\ Z(\chi)/Z(G)=p$. This gives that either $Z(\chi)=Z(G)$ or $|Z(\chi)/ Z(G)|=p$. If $Z(\chi)=Z(G)$ we get that $|G/Z(\chi)|=|G/Z(G)|=p^{2n+1}$, which contradicts that $\chi(1)^{2}=|G/Z(\chi)|$. Hence $|Z(\chi)/Z(G)|=p$.\\

  $(iii)\Rightarrow (ii)$ Obvious.\\

   $(iii)\Rightarrow (i)$  By (ii) We get that the group $G$ satisfies the weak condition on normal subgroups of $G$ and hence by Theorem \ref{thm1} we get that exp $G'=p$ or $p^{2}$. If exp $G'=p^{2}$ it follows from Theorem \ref{thm1} that $|G/Z(G)|=p^{4}$. This contradiction gives that exp $G'=p$. Since $|Z(\chi)/Z(G)|=p$, we get that $Z(\chi)/Z(G)$ is cyclic.

\end{proof}

\end{thm}

In Theorem C of \cite{fernandez2001groups} the authors found that for any $p$-group $G$ satisfying the weak condition on normal subgroups, $cd(G)=\{1, p^{n}\}$ if $|G/Z(G)|=p^{2n+1}$. In the following theorem we obtain information on character degree sets of such groups under the condition that  $|G/Z(G)|=p^{2n}$.

\begin{thm}

Let $G$ be a finite $p$-group of nilpotency class not equal to $3$  and $|G/Z(G)|=p^{2n}$. If $G$ satisfies the weak condition on normal subgroups then $cd(G) \subseteq \{1, p^{n-1}, p^{n}\}$.

\begin{proof}
By Lemma \ref{lemma4} we get that $c(G)\leq 4$.\\
First we consider the case that the nilpotency class of $G$ is $2$.
Let $K$ be the kernel of a non-linear irreducible character of $G$. First we assume that $|KZ(G)/Z(G)|= p$ . Since $KZ(G)/Z(G)\leq Z(\chi)/Z(G)$, by Lemma \ref{lemma3.1} we get that the Abelian group $Z(\chi)/Z(G)$ contains a subgroup isomorphic to $Z(\chi)/KZ(G)$, which is cyclic. Since $|KZ(G)/Z(G)|= p$, either $|Z(\chi)/KZ(G)|=1$  or $Z(\chi)/KZ(G)$  is of maximal order. If $|Z(\chi)/KZ(G)|=1$, then we get $Z(\chi)=Z(G)$ and hence $|Z(\chi)/Z(G)|=p$. Consequently it follows that $\chi(1)^{2}=|G/Z(\chi)|=|G/Z(G)|=p^{2n}$.  If $|Z(\chi)/KZ(G)|$  is of maximal order then by the proof of 2.1.2 of \cite{kurzweil2004finite} we get that $Z(\chi)/Z(G)\cong C_{p^{k}}\times C_{p}$, where $p^{k}=|Z(\chi)/KZ(G)|$. Since $G$ satisfies the weak condition on normal subgroups,   by Theorem \ref{thm1} we get that exp $\ G'=exp\ G/Z(G)=p$ or $p^{2}$. This gives that exp $\ Z(\chi)/Z(G)=p$ or $p^{2}$. If exp $\ Z(\chi)/Z(G)=p$ we get that $k=0$ or $k=1$. If $k=1$, then $|Z(\chi)/Z(G)|=p^{2}$ and hence $\chi(1)^{2}=|G/Z(\chi)|=p^{2n-2}$. If $k=0$ we get that 
  $Z(\chi)/Z(G)\cong  C_{p}$. Since $|KZ(G)/Z(G)|=p$ and $KZ(G)/Z(G) \leq Z(\chi)/Z(G)$, we get that $KZ(G)/Z(G)=Z(\chi)/Z(G)$ and this in turn gives that $|Z(\chi)/Z(G)|=p$. Thus $|G/Z(\chi)|=p^{2n-1}$ and this contradicts that  $\chi(1)^{2}=|G/Z(\chi)|$.
  If exp $\ G'=p^{2}$ we get that $Z(\chi)/Z(G)\cong C_{p^{2}}\times C_{p}$. This gives that $|Z(\chi)/Z(G)|=p^{3}$ which in turn gives that  $\chi(1)^{2}=|G/Z(\chi)|=p^{2n-3}$, a contradiction.
  
  Next assume that $|KZ(G)/Z(G)|= 1$ so that $K\leq Z(G)$. Then $Z(\chi)/Z(G)$ is cyclic. Since $exp\ Z(\chi)/Z(G)=p$ or $p^{2}$, it follows that $Z(\chi)=Z(G)$ or $Z(\chi)/ Z(G)\cong C_{p}$ or $C_{p^{2}}$. If $Z(\chi)/ Z(G)\cong C_{p}$ we get that $\chi(1)^{2}=|G/Z(\chi)|=p^{2n-1}$, a contradiction. Hence $Z(\chi)=Z(G)$ or $Z(\chi)/ Z(G)\cong C_{p^{2}}$. Thus $\chi(1)^{2}=|G/Z(\chi)|=p^{2n}$ or $p^{2n-2}$.\\
  Next we consider the case that the nilpotency class of $G$ is $4$. By Theorem F of \cite{fernandez2001groups} we get that $|G/Z(G)|=p^{4}$ and hence $cd(G)\subseteq \{1, p, p^{2}\}$.

\end{proof}

\end{thm}

\end{document}